\documentclass[14pt]{amsart}

\usepackage{amssymb, array, verbatim, amscd}
\usepackage{amsmath, amsfonts,amsthm}

\theoremstyle{plain}
\newtheorem{theorem}{Theorem}[section]
\newtheorem{corollary}[theorem]{Corollary}
\newtheorem{lemma}[theorem]{Lemma}
\newtheorem{proposition}[theorem]{Proposition}

\theoremstyle{definition}
\newtheorem{algorithm}[theorem]{Algorithm}
\newtheorem{definition}[theorem]{Definition}
\newtheorem{example}[theorem]{Example}

\newtheorem{remark}[theorem]{Remark}

\newcounter{abc}

\newcounter{ABC}

\newcommand{\Hom}{\operatorname{Hom}}

\newcommand{\la} {\langle}
\newcommand{\ra} {\rangle}
\newcommand{\N}{{\mathbb N}}
\newcommand{\Z}{{\mathbb Z}}

\newcommand{\R}{{\mathcal R}}
\newcommand{\E}{{\mathcal E}}

\newcommand{\J}{{\mathcal J}}

\newcommand{\I}{{\mathcal I}}

\begin{document}
\title[A new algorithm for finding the nilpotency class of a finite $p$-group] {A new algorithm for finding the nilpotency class of a finite $p$-group
describing the upper central series}

\author{Maria A. Avi\~no-Diaz}
\address{Department of Mathematic-Physics\\ University of Puerto Rico\\
Cayey, Puerto Rico 00736 } \email{mavino@uprr.pr}
\date{\today}
\begin{abstract}
 In this paper we describe an algorithm for finding
the nilpotency class, and the upper central series of the maximal
normal $p$-subgroup $\Delta (G)$ of the automorphism group, $Aut(G)$
of a bounded (or finite) abelian $p$-group $G$. This is the first
part of two papers devoted to compute the nilpotency class of
$\Delta(G)$ using formulas, and algorithms that work in almost all
groups. Here, we prove that for $p\geq 3$  the algorithm always
runs. The algorithm describes a sequence of ideals of the Jacobson
radical, $\J$, and because $\Delta(G)=J+1$, this sequence induces
the upper central series in $\Delta (G)$.
\end{abstract}

\maketitle

  \section{Introduction}
The automorphism group of an abelian $p$-group, was studied for  K.
Shoda in 1928, under the advice of Emmy Noether, and he gave the
description of the endomorphism ring and a characterization of the
automorphism group using a matrix representation  over the integer
modulo a primary number $\Z_{p^n}$, \cite{S28}.  The maximal normal
$p$-subgroup of the automorphism group of a bounded  abelian
$p$-group $G$ plays a very important role in the description of  the
automorphism group, because $Aut (G)$ decomposes on semidirect
product of $\Delta (G)$ in several cases, see \cite{A4}.  This paper
is the first part of two papers devoted to the description of the
nilpotency class of the maximal normal $p$-subgroup $\Delta (G)$ of
the automorphism group of a bounded abelian $p$ group
$G=\oplus_{i\leq s} G_i$, where $G_i$ is a homocyclic group of
exponent $p^{n_i}$ and  p-rank $r_i$, with $n_i<n_{i+1}$ for all
$i$. Here, we prove that some special sequence of ideals, the upper
annihilating series of the Jacobson radical $\mathcal J$ of the
endomorphisms ring of $G$ plays a very important role in the
description of the upper central series of the maximal normal
p-subgroup $\Delta (G) =\mathcal J +1 \leq Aut (G)$. In fact, we
obtain the nilpotency class  of $\Delta (G)$, and characterize its
upper central series (ucs) $\{\mathcal Z_t\}$. to do this, we use
the upper annihilating sequence (uas) $\{J_t\}$, where
$\J_t$=Annihilator of $(\J/\J_{t-1})\}$, for all $t$, a sequence of
two sided ideals of $\J$  in the endomorphism ring $\mathcal E (G)$.
This sequence was defined and studied for finite  groups in \cite
{A1, A2,A3, A0}, and for bounded groups in \cite{AS01}. Here we
obtain new results for the (uas), in Section \ref{(uas)}.
 For most the cases this sequence determine the upper central series.
 The nilpotency class of a finite $p$-group is usually computing the lower central series,  but here we introduce
 an algorithm for computing the upper central series and of course the nilpotency class.
 We have a recursive function that
permits the construction of the (ucs). One of the most interesting
results is that the description, and the length of the upper central
series only depends on the exponents of the group $G$, and the rank
$r_s$ of the homocyclic subgroup $G_s\cong
\left(\Z_{p^{n_s}}\right)^{r_s}$ of maximal exponent in the group
$G$. This is true even if the group is infinite but bounded. Here,
we introduce a method to construct central series associated to
annihilating series, that is, if we have a ring $R$ with identity
$1$, and there exists a bilateral ideal $I$  that has a finite
annihilating sequence then $H=I+1$ is a nilpotent group with degree
less or equal the length of the upper annihilating sequence of $I$.

A group $G$ has the (uas)-(ucs)-property for the Jacobson radical
$\J$ of $\E$, if the upper central series $(\mathcal Z_t)_{t\geq 1}$
of $\Delta(G)$ satisfies the following $\mathcal Z _t=(\J_t+1)
Z_\Delta$ for all $t$ from $1$ to $n(\Delta)$=the nilpotency class
of $\Delta(G)$, and $ Z_\Delta $ is the intersection of the subgroup
$\Delta (G)$ with the center of $\E$.

  Associated to the Upper Annihilating
Sequence  of $\J$, we have a special function that is defined in a
recursive form, because the function calls itself to construct the
ideals, we have the initial value, and the condition to stop given
by the the ideal function associated to the radical $\J$, and the
Theorem \ref{stop}. We call this function  Upper  Function and we
use it to describe an algorithm for finding the Upper Central Series
of $\Delta (G)$, it is done in  Section \ref{Algorithm1}. In this
paper, the following theorem is proved:
\begin{theorem} Let $G\cong \oplus_{i=1}^s G_i$ be a bounded abelian $p$-group,
 where the $p$-rank of $G_i\cong \left(\Z_{p^{n_i}}\right)^{r_i}$
is the ordinal number  $r_i$,  where $0<n_1<\cdots < n_s$. Let
$\{\J_t\}$ be the upper annihilating series of the radical $\J$ of
$\E$, and let $\{\mathcal Z_t\}$ be the upper central series of the
maximal normal $p$-subgroup of $Aut (G)$, $\Delta (G)$. If
$Z=\mathcal Z \cap \Delta (G)$ and one of the following cases holds:
\begin{itemize}
\item[(1)]  $p\geq 3$
\item [(2)]$p=2$, $\sigma (G)\geq 2$,
\item [(3)] $p=2$, $r_s>1$, and $s\geq 2$.
\end{itemize}
then  for $t\geq 1$, $\begin{array}{|c|} \hline \mathcal
Z_t=(\J_t+1)Z.\\  \hline\end{array}$  The description of the
elements in the hypercenters $\mathcal Z _t$ is the following:
\[ \mathcal Z _t=\left\{ C+1_{rr}\in \Delta (G)| \left\{\begin{array}{ll}
c_{ij}^{(l,k)}\equiv 0 \pmod{p^{n_j-\alpha{(i,j,t)}}}& i\ne j;\\
c_{ii}^{(l,l)}\equiv c_{ss}^{(r,r)} \pmod {p^{n_i-\alpha{(i,i,t)}}}& for \ i\leq s \\
c_{ii}^{(l,k)}\equiv 0 \pmod{p^{n_i-\alpha{(i,i,t)}}}& l\ne k\cr
\end{array} \right\}\right\}  \] where $C=(c_{ij}^{(l,k)})_{r\times r}$ is a matrix with
entries by columns in $\Z_{p^{n_j}}$, where $r=\sum_{i=1}^sr_i$.
\end{theorem}

In the second paper,  general formulas for the length of the upper
annihilating sequence, and the nilpotency class of $\Delta$ are
given.
 \section{Preliminaries and Notation}\label{Prel}
 Let $G=\oplus_{i=1}^s G_i$ be a decomposition of a bounded abelian
$p$-group $G$ into homocyclic subgroups $G_i\cong
\left(\Z_{p^{n_i}}\right)^{r_i}$, where $0<n_1<\cdots < n_s$ are
integers and $r_i$ cardinals.

 The endomorphism ring $\E$ of $G$ will be represented as an
$s\times s$ matrix ring $({\E}_{ij})$, where for all $(i,j)$,
${\E}_{ij}=\Hom(G_i,G_j)$, considered as an
${\E_i}$-${\E_j}$-bimodule.
   For each endomorphism  $f=(f_{ij})_{s\times s}\in \E$ we consider the
 functions $f_{ij}$ defined as follows. If $x\in G$ decomposes
   as sum of elements $x_i\in G_i$, then $f(x)=\sum_i f(x_i)=\sum_i\sum_j x_{ij}$,
   where $f(x_i)=\sum_j x_{ij}\in G_j$, so $f_{ij}(x_i)=x_{ij}$, and
   $f_{ij}\in \E_{ij}$.
   So, $\E\cong(\E_{ij})$. By the above decomposition we have that
   the matrix $\E_{ij} $ has entries in $\Z_{p^j}$ and satisfies the following
condition $\E_{ij}\equiv 0 \pmod{p^{n_j-n_i}}$ if  $i<j$.
 The center $\mathcal Z $ of ${\E}$ is the ring of scalar matrices
$c\mathcal {I}_\E$ where $ c_{ii}\equiv c \pmod {p^{n_i}}$,
 and for $i<j,\ c_{ii}\equiv c_{jj}\pmod{p^{n_j}}$. The maximal normal
  $p$-subgroup of $Aut (G)$ is denoted by $\Delta(G)=1+\J$.

 Here, we always consider an ideal like a bilateral ideal. $I$ is an ideal of ${\E}$ if and only if $I=(I_{ij})$ where each
$I_{ij}$ is a sub-bimodule of ${\E}_{ij}$.
 $I$ is an ideal of ${\E}$ if and only if there exists an
ideal function $\beta$ such that the matrix representation
 of the ideal is the following $I=(p^{\beta(i,j)}E_{ij})$, where $p^{\beta(i,j)}E_{ij}=I_{ij}$, and
$E_{ij}$ is a $r_i\times r_j$-matrix with entries in the integers
modulo $p^{n_j}$,
 \cite{A0,AS01}.

${\E}$ has Jacobson radical $\J=(\J_{ij})$ where for all $i\ne j,\
\J_{ij}={\E}_{ij}$ and $\J_{ii}= p{\E}_{ii}$, \cite{F,J}. We denote
the ideal function associated to the Jacobson Radical by
$n_j-f_\J(i,j)$, where $f_\J(i,j)=\left\{\begin{array}{ll}
n_{i\wedge j} & if \
i\ne j\\
n_j-1 & if\ i=j\\
\end{array}\right. ,$
$min\{i,j\}=i\wedge j$.

 We will say that $G$ has  underlying type
$\underline{t}(G)=(n_1,\dots, n_s)$, with $n_i< n_{i+1}$, for all
$i$. We denote by $\sigma(G)$ the minimum gap in the sequence
$n_1<n_2<\cdots <n_s$, that is $\sigma(G)=min_j(n_j-n_{j-1})$
 Let $\ell=\max \{j:n_j-n_{j-1}=\sigma(G)\}$.

If the rank of $G$ is denoted by $r=\sum_i r_i$, then let
$G=\oplus_{i=1}^n\oplus_{k=1}^{ r_i}\la x_{ik}\ra$ be a fixed
decomposition of $G$ into cyclic summands, and let $X=\cup X_i$
where $X_i=\{x_{1i},\ldots , x_{r_ii}\}$ is a basis of $G_i$.   We
define a set of endomorphisms in $\J$, called \textbf{elementary
endomorphisms}, that are the common elementary matrices, but
consider in the radical $\J$,
\begin{itemize}
\item[(1)] For all $x_{ki},\, x_{li}\in X_i,\ e^{(k,l)}_{ii}$
maps $x_{ki}$ to $px_{li}$ and annihilates the complement of $\la
x_{ki}\ra$ with respect to the basis $X$. The elementary matrix of
$e^{(k,l)}_{ii}$ has $p$ in the place $kl$ of the diagonal cell
$ii$. For example, the endomorphism $e^{({k,l})}_{11}$ has matrix
representation
\[ E_{1,1}^{(k,l)}=\left( \begin{array}{c||c|c}
pE^{kl}_{11}& 0& 0 \\
\hline \hline
0 &0 & 0 \\
\hline
0&0&0  \\
\end{array}\right),\]
where the matrix $E^{kl}_{11}$ is the elementary $r_1\times
r_1$-matrix which has the number $1$ in the place $kl$.
\item[(2)] For $i>j$, and for all $x_{ki}\in X_i,\ x_{lj}\in X_j$, $e^{(k,l)}_{ij}$
maps $x_{ki}$ onto $x_{lj}$ and annihilates the complement of $\la
x_{ki}\ra$. The elementary matrix of $e^{(i,l)}_{ij}$ has $1$ in the
place $kl$ of the  cell $ij$. For example, the endomorphism
$e^{(kl)}_{21}$ has matrix representation
\[ E_{2,1}^{(k,l)}\left( \begin{array}{c||c|c}
0& 0& 0 \\
\hline \hline
E^{kl}_{21}&0 & 0 \\
\hline
0&0&0  \\
\end{array}\right),\]
where the matrix $E^{kl}_{21}$ is the elementary $r_2\times
r_1$-matrix which has the number $1$ in the place $kl$.
\item[(3)] For $i<j$, and for all $x_{ki}\in X_i,\ x_{lj}\in X_j$ with $i<j$, $e^{(k,l)}_{ij}$ maps
$x_{ki}$ onto $p^{n_j-n_i}x_{lj}$ and annihilates the complement of
$\la x_{ki}\ra$. The elementary matrix of $e^{(k,l)}_{ij}$ has
$p^{n_j-n_i}$ in the place $kl$ of the  cell $ij$. For example, the
endomorphism $e^{(k,l)}_{12}$ has matrix representation
\[ E_{1,2}^{(l,k)}=\left( \begin{array}{c||c||c}
0& p^{n_2-n_1}E^{kl}_{12}& 0 \\
\hline \hline
0&0 & 0 \\
\hline
0&0&0  \\
\end{array}\right),\]
where the matrix $E^{kl}_{12}$ is the elementary $r_1\times
r_2$-matrix which has the number $1$ in the place $kl$.
\end{itemize}
\section{Central Series induced by  Annihilating
Sequences }\label {(uas)}
  Let $I$ be an ideal of a ring  $R$. Let
    $A_1=Ann (I)=\{a\in I : aI=0=Ia\}$ be  the annihilator of $I$.
  The following definitions appear in \cite{A1,A2,A3,AS01}.
\begin{definition} The upper annihilating sequence, (uas), of the ideal $I$ is defined by \[ I_0
=0,\ and \ for\  t\geq 1,\ I _t /I_{t-1}=Ann (I /I_{t-1}).\] In
other words, $I_t=\{a\in I : aI\subseteq I_{t-1}$, $I a\subseteq
I_{t-1}\}$. This is  an ascending sequence of ideals $I_t$. If there
is a positive integer $d$ such that
  \[0=I_0\subset  \cdots  \subset  I_{t-1} \subset
U_{t } \subset  \cdots  \subset  I_d= I,\]
   we say that  the annihilating length of $I$ is $d$.
\end{definition}
   We  make a generalization of the above definition in order to find the
   ideal function associated to the upper annihilating sequence of the Jacobson radical of $\E$,
   This definition will be used in the next section in a ring R with unit 1.
\begin{definition}
   We will say that a sequence of ideals $A_t$ is an annihilating
   sequence of $I$ if \[ A_0 =0,\  A_1=Ann(I)\ and \  for \
   t\geq 2,\  A _t /A_{t-1}\subseteq Ann (I /A_{t-1}),\] that is
   $A_t=\{ a\in I |ab \in A_{t-1},ba \in A_{t-1},\forall b\in I\}.$
   So, \[0=A_0\subset  A_1=Ann(I)\subset \cdots  \subset  A_{t-1} \subset A_{t } \subset  \cdots  \subset  A_{l(A_t)}= I.\]
   If for some $t$ the ideal $I=A_t$, then the first $t$ with this property
   will be called the length of $A_t$  and  denoted by
   $l(A_t)$, \cite{A0}.
   \end{definition}
It is clear that $l(I_t)\leq l(A_t)$ for all annihilating sequence
$(A_t)_t$ of $I$.
 \subsection*{Central series induced by annihilating
sequences} In this section we introduce the mean relation between
the upper annihilating series of the radical $\J$, and the central
series of the maximal normal $p$-subgroup $\Delta (G)=\J +1$. This
relation is given in the Theorem \ref {uas=ucs}.

The length of the (uas) of $I$ is called the annihilating length of
$I$.

If $R$ is a ring with 1, and $I$ is an ideal of $R$, with finite
annihilating length, then $1+I=\Gamma$ is a normal subgroup  of the
group of units of $R$, that satisfies the following theorem.
\begin{theorem}\label{uas=ucs} Let $\{A_t\}$ be an annihilating sequence of
the ideal $I$ of length $l(A_t)$, and let $\Gamma=I+1$. If $Z=
\mathcal Z \cap\Gamma$, where $\mathcal Z$ is the center of $R$.
Then
\[\{1\}\leq (A_1+1)Z\leq \cdots\leq
(A_t+1)Z\leq \cdots\leq  (A_{l(A_t)}+1)Z=\Gamma \] is a central
series for $\Gamma$, and $\Gamma $ is a nilpotent group.
Consequently, the nilpotency class of $\Gamma$ is less than or equal
to the annihilating length of $I$.
\end{theorem}
\begin{proof} Let  $\Gamma _0={1}$, and $\Gamma _t=(A_t+1)Z$ for $1\leq t\leq l(A_t)$.
If $\alpha=a+1\in\Gamma_t$ and $\beta=b+1\in\Gamma=I+1$, then we
will prove that the commutator  $[\alpha,\beta]\in\Gamma_{t-1}$.
Computing $[\alpha,\beta]$ we have that
\[\begin{array}{lll}
   [\alpha,\beta]&=& \alpha^{-1}\beta^{-1}\alpha \beta=
\alpha^{-1}\beta^{-1}(\alpha
\beta-\beta\alpha)+1=\\
&=&\alpha^{-1}\beta^{-1}((a+1)(b+1)-(b+1)(a+1))+1=\alpha^{-1}\beta^{-1}(ab-ba)+1
\end{array}\]
Denoting $\alpha=(1+a_t)z$ with $a_t\in A_t$ and $z\in Z$, we have
that $\alpha=\overline {a}+z$, with   $\overline{a}\in A_t$. On the
other hand \[ab-ba=(\overline {a}+z)(b+1)-(b+1)(\overline
{a}+z)=\overline {a}b-b\overline {a}=(a,b)\in A_{t-1},\] because
$\overline{a}\in A_t$. So $\alpha^{-1}\beta^{-1}(ab-ba)\in A_{t-1}$
and $[\alpha,\beta]\in A_{t-1}+1\subseteq \Gamma _{t-1}$.
\end{proof}
\begin{corollary}\label{COR}
 Let $\{A_t\}$ be an annihilating sequence of
the ideal $I$ of $R$ and let $\Gamma=I+1$. If $Z=\mathcal Z
\cap\Gamma$, where $\mathcal Z$ is the center of $R$. Then the
central series $(\Gamma _t)_{t\ge 1}$, and $\Gamma _0=\{1\}$,
satisfies that the following property, for all $t>0$:
\[If \ a+1\in \Gamma _t=(A_t+1)Z,\ then \ (a,b)\in A_{t-1}\ for \
all \  b+1\in \Gamma\]
\end{corollary}

Of course, these results are true for the endomorphism ring $\E,$
for the radical $J$, and the maximal normal p-subgroup of $Aut G$,
$\Delta (G)= \J+1$.
\section{The annihilating
functions} In this section we introduce the function, that we call
Upper annihilating  function which describes the upper annihilating
sequence. In particular the function \ref{anni} gives the algorithm
to construct the upper central series of $\Delta (G)$. First we
compute the annihilator of $\J$, $Ann(\J)=\J_1$.
 \begin{lemma}\label{anni1}
The annihilator of $\mathcal J$ is the ideal $\J_1$ described by the
matrices  $(A_{ij})_{s\times s}$, such that $A_{ij}=0$ for all
$(i,j)\ne (s,s)$, and $(A_{(s,s)})=p^{n_s-1}\E_{ss}$.
\end{lemma}
\begin{proof}By definition $\J_1=Ann(\J)=\{A\in \J|AB=BA=0, \forall B\in
\J\}$ using the matrix representation of $\E$. Suppose that $A\in
\J_1$, and $B\in \J$, then the condition for the annihilating ideal
implies that
\begin{equation}\label{J-1}(i,j)\ \hspace{1.2in}
\sum_{k=1}^s A_{ik}B_{kj}\equiv \sum_{k=1}^sB_{ik}A_{kj}\equiv 0
\pmod {p^{n_j}}.
\end{equation}
If  $B_{ks}=0$, for $k\ne s$, and $B_{ss}=p\I_{ss}$, then for all
$i\ne s$ in $(\ref{J-1}-(i,s))$ , we have that $A_{is}\equiv 0
\pmod{p^{n_s}}$,  and $pA_{ss}\equiv 0 \pmod {p^{n_s}}$, so
$A_{ss}\equiv 0 \pmod {p^{n_s-1}}$. Similarly we prove that
$A_{ij}\equiv 0 \pmod {p^{n_i}}$, for all $(i,j)\ne (s,s)$. On the
other hand, an element $A\in\J$ with these conditions satisfies that
$A\in \J_1$. Then our claim holds. We can observe that the ideal
function associated to $\J_1=Ann(\J)$ is $\beta (s,s)=n_s-1$, and
$\beta (i,j)=0$, for $(i,j)\ne (s,s)$.
\end{proof}
We use  the matrix representation of $\mathcal J=({p^{n_j-n_{i\wedge
j}+\delta _{ij}} \mathcal E }_{ij})$, where $\delta _{ij}$ is the
Kronecker's delta, and $ \mathcal E_{ij}$ is the set of $r_i\times
r_j$-matrices with entries in $\Z_{p^{n_j}}$.
\begin{theorem}\label{anni}
The sequence of ideals
\[0=A_0\subset  \cdots  \subset  A_{t-1} \subset
A_{t } \subset  \cdots  \subset  A_n= \mathcal J,\] is an
annihilating sequence  of $\mathcal J$ if, and only if, there exists
an ideal function $f$ such that $ A_t=(p^{n_i-f{(i,j,t)}}\mathcal
E_{ij}),$ for all $(i,j,t)\in D$, which satisfies:
\begin{enumerate}
\item $f(s,s,1)=1$, and $f(i,j,1)=0$,
\item $f(i,j,t)\leq f(i,j,t-1)+1,$
\item $f(i-1,j,t)\leq f(i,j,t-1)\leq f(i,j,t),$
\item $f(i,j-1,t)\leq f(i,j,t-1)\leq f(i,j,t),$
\item $f(i+1,j,t)\leq n_{i+1}-n_i+ f(i,j,t-1),$
\item $f(i,j+1,t)\leq n_{j+1}-n_j+ f(i,j,t-1).$
\end{enumerate}
We will call this class of functions, the annihilating functions.
\end{theorem}
\begin{proof} ($\Rightarrow$ ) Suppose
$\{A_t\}$ is an annihilating sequence of $\mathcal J$. So, there
exists an ideal function associated to this sequence that we denote,
$\beta_t(i,j)=n_i-f(i,j,t)$. We claim that $f$ is an annihilating
function, so we prove the six properties of the function $f$ in the
theorem. By Lemma \ref{anni1} the condition (1) holds. Because
$A_{t-1}\subseteq A_t$, for $t\geq 2$, we have that
\[p^{n_i-f{(i,j,t-1)}}\mathcal E _{ij}\subseteq p^{n_i-f{(i,j,t)}}\mathcal E _{ij} \Rightarrow f{(i,j,t-1)}\leq  f{(i,j,t)}. \]
By definition of annihilating sequence,  if we take $A\in A_t$, and
$B \in \mathcal J$, then $AB\in A_{t-1}$. On the other hand,
$AB=(\sum_{k=1}^s A_{ik}B_{kj})_{s\times s}$, with $A_{ik}\in
p^{n_k-f(i,k,t)}\mathcal E_{ik}$, and $B_{kj}\in p^{n_j-n_{k\wedge
j}+\delta _{kj}} \mathcal E _{kj}$. For all $1\leq k\leq s$,we have
the following
\[ A_{ik}B_{kj}\in p^{n_k-f(i,k,t)+n_j-n_{k\wedge j}+\delta
_{kj}}\mathcal E_{ik}\mathcal E_{kj}\Rightarrow\]
\[p^{n_k-f(i,k,t)+n_j-n_{k\wedge j}+\delta _{kj}}\mathcal E_{ij}\subseteq
p^{n_j-f(i,j,t-1)}\mathcal E_{ij}\hbox{ because } AB \in A_{t-1}.
\]
Therefore
\[\  n_k-f(i,k,t)+n_j-n_{k\wedge j}+\delta
_{kj}\geq n_j-f(i,j,t-1)\Rightarrow \]
\begin{equation}\label{eq1}f(i,k,t)-\delta _{kj} \leq n_k-n_{k \wedge
j}+f(i,j,t-1). \end{equation} On the other hand, $BA\in A_{t-1}$,
and   $BA=(\sum_{k=1}^s B_{ik}A_{kj})_{s\times s}$, where the matrix
$B_{ik}\in p^{n_k-n_{i\wedge k}+\delta _{ik}}\mathcal E_{ik}$, and
$A_{kj}\in p^{n_j-f(k,j,t)}  \mathcal E_{kj}$. For all $1\leq k\leq
s$, we have the following
\[\  n_j-f(k,j,t)+n_k-n_{i\wedge k}+\delta _{ik}\geq n_j-f(i,j,t-1)\Rightarrow \]
\begin{equation}\label{eq2}f(k,j,t)-\delta _{ik} \leq n_k-n_{i \wedge k}+f(i,j,t-1).
\end{equation}
Using the equation (\ref{eq1}) we obtain the following relations
\[\begin{array}{|r|c|}
          \hline
k=j+1 \Rightarrow & f(i,j+1,t)\leq  n_{j+1}-n_{j}+f(i,j,t-1) \hspace{.4 in}   (6)  \\
\hline k=j\Rightarrow & f(i,j,t)\leq f(i,j,t-1)+1 \hspace{1.2
in}(2)\\\hline
 k=j-1 \Rightarrow & f(i,j-1,t) \leq (i,j,t-1) \hspace{1.3in} (4)\\\hline
\end{array}\]
Using the equation (\ref{eq2}) we obtain the following relations.
\[\begin{tabular}{|r|c|}
          \hline
$k=i+1 \Rightarrow $& $f(i+1,j,t)\leq n_{i+1}-n_i+ f(i,j,t-1)$ \hspace{.3 in} (5)  \\
\hline $k=i\Rightarrow$& $f(i,j,t)\leq f(i,j,t-1)+1$
\hspace{1.1in}(2)\\\hline
 $k=i-1 \Rightarrow$ &$ f(i-1,j,t) \leq f(i,j,t-1)$\hspace{1.1in} (3)\\\hline
\end{tabular}\]
The condition(1) is proved by Lemma \ref{anni1}. So, our claim
holds.

($\Leftarrow$) It is easy to prove that the chain of ideals
$A_t=(p^{n_i-f(i,j,t)}\mathcal E_{ij})$, where $f$ is an
annihilating function, is an annihilating sequence of $\mathcal J$.
So, the theorem holds.
\end{proof}
It is clear that the upper annihilating sequence  has the minimal
length  between the annihilating sequences. Therefore, the
annihilating function associated to the (uas) is the minimum
annihilating function, that we denote by $\alpha $. The symbol
$\bigwedge $ means the minimum number.
\begin{corollary}[The Upper Function]\label{formulas}
For $t\geq 1$, the upper annihilating function $\alpha$ is a
recursive function on $t$,  defined by the  following formulas:
\begin{itemize}
\item[(1)] $\alpha(s,s,1)=1,$  $\alpha(i,j,1)=0,$  for $(i,j)\ne (s,s)$
\item[(2)] $\alpha{(i,j,t)}=\bigwedge
\{n_j-n_{j-1}+\alpha{(i,j-1,t-1)}; n_i-n_{i-1}+\alpha{(i-1,j,t-1)};
$\\
$\alpha{(i,j,t-1)}+1;\alpha{(i+1,j,t-1)};\alpha{(i,j+1,t-1)}\}.$
\end{itemize}
we consider the function  $\alpha:S\times S\times T \rightarrow \N$,
where $S=\{1,\ldots,s\}$, and $T=\{1,\ldots ,d=l(\J_t)\}$.
\end{corollary}
\begin{proposition}\label{prop:uf}
 The following properties holds
\begin{itemize}
\item [(1)] $\alpha(i+1,j,t)\leq \alpha (i,j,t)$, for all $1\leq i\leq
s-1$,
\item[(2)] $\alpha(i,j,t)\leq \alpha (i,j,t+1)$, for all $t\leq
l(\J_t$,
\item[(3)]$\alpha(i,j+1,t)\leq \alpha (i,j,t)$, for all $1\leq j\leq
s-1$,
\item[(4)]  For all $j<s$, $n_j-\alpha(j,j,t)\leq
n_s-\alpha(s,s,t)$.
\end{itemize}
\end{proposition}
\begin{proof}
Because $\alpha$ is a annihilating function, it satisfies the
properties in Theorem \ref{anni}, so again the properties (1), (2),
(3).

The proof of property (4) is the following
\[\alpha(j+1,j+1,t)\leq \alpha(j,j+1,t)\leq \alpha(j,j,t)\]
by properties (1), and (3). On the other hand, we have that
\[n_{j+1}-n_j\geq 0, \ and \ \alpha(j,j,t)-\alpha(j+1,j+1,t)\geq
0.\] Therefore, $n_j-\alpha(j,j,t)\leq n_{j+1}-\alpha(j+1,j+1,t)$,
for all $j< s$, and the property (4) holds.

\end{proof}
\begin{remark}
 Using the matrix
representation of $\mathcal E$, we proved that
$J_t=(p^{n_i-\alpha{(i,j,t)}}\mathcal E_{ij})$,we want to remark
that the integers $\alpha{(i,j,t)}$ depend only on the exponents
$(n_1,\dots, n_s)$ of the group $G$, but not the ranks of the
homocyclic components of $G$.
\end{remark}
In this paper we use the following notation, for $x\in \R_{\geq 0}$,
$[x]=$the integer part of $x$, if $x<0$, then $[x]$ means 0.
\begin{lemma}[\textbf{Case 1 for $\alpha$}] If $\sigma(G)\geq 2$ then
$\alpha{(i,j,t)}=[t+i+j-2s].$
\end{lemma}
\begin{proof}
If  $\sigma(G)\geq 2$, we use induction in order to prove that
$\alpha{(i,j,t)}=[t+i+j-2s].$ For $t=1$ is trivial, suppose that the
theorem holds for $t-1$, then
\[\alpha{(i,j,t)}=\bigwedge \{n_j-n_{j-1}+[t-2+j+i-2s];n_i-n_{i-1}+[t-2+j+i-2s];\]
\[[t-1+i+j-2s]+1;[t+i+j-2s];[t+i+j-2s]\}=[t+i+j-2s]\]
\end{proof}
\begin{lemma}[\textbf{Case 2 for $\alpha$}]\label{alpha}  If $\sigma (G)=1$, and $n_s-n_{s-1}=1$ then
\[\alpha{(i,j,t)}=\left[\displaystyle{\frac{t+i+j-2s+1}{2}}\right].\]
\end{lemma}\begin{proof}
We use induction to prove the Lemma, so for $t=1$ is trivial.
Suppose the lemma holds for $t-1$, and we have
\[\alpha{(i,j,t)}=\bigwedge\left\{n_j-n_{j-1}+\left[\displaystyle{\frac{t+i+j-2s-1}{2}}\right]
;n_i-n_{i-1}+\left[\displaystyle{\frac{t+i+j-2s-1}{2}}\right];
\right. \] \[\left.
\left[\displaystyle{\frac{t+i+j-2s}{2}}\right]+1;\left[\displaystyle{\frac{t+i+j-2s+1}{2}}\right];
\left[\displaystyle{\frac{t+i+j-2s+1}{2}}\right] \right\}\] It is
obvious that
$\alpha(i,j,t)=\left[\displaystyle{\frac{t+i+j-2s+1}{2}}\right]$
\end{proof}
The third case for $\alpha$ is the following: $\sigma (G)=1$, and
$n_j-n_{j-1}=1$, but $j<s$. There are formulas for all the cases. In
the second paper,  general formulas for the length of the upper
annihilating sequence, and the nilpotency class of $\Delta$ are
given.
\begin{example}\label{Ex1} Consider a group of type $\underline t(G)=(3,5,6,8,10)$
with $p$-rank $r_i=1$, for all $i$. The matrices in the annihilating
sequence have the entries by columns in $\Z_{p^3}$, $\Z_{p^5}$,
$\Z_{p^6}$, $\Z_{p^8}$, and $\Z_{p^{10}}$ respectively. Observe that
this is the third  case, because $\sigma (G)=1$, but
$n_s-n_{s-1}=2$. Using the Corollary \ref{formulas}, the upper
annihilating sequence is the following, we compute the recursive
formula for the ideals $\J_t$. We only include the correspondent
power of $p$ in the place$(i,j)$, the meaning is: that place  is
congruent to  $0$ modulo $p^x$ . The matrices represent the sequence
of ideals $\J_1\subset \cdots \subset \J_{14}=\J$.
\[\left[\begin{array}{ccc|cc}
0&0&0&0&0\\
0&0&0&0&0\\
0&0&0&0&0\\
\hline
0&0&0&0&0\\
 0&0&0&0&p^9\cr
\end{array}\right], \
\left[\begin{array}{ccc|cc}
0&0&0&0&0\\
0&0&0&0&0\\
0&0&0&0&0\\
\hline
0&0&0&0&p^9\\
0&0&0&p^7&p^8\cr
\end{array}\right], \left[\begin{array}{ccc|cc}
0&0&0&0&0\\
0&0&0&0&0\\
0&0&0&0&p^9\\
\hline
0&0&0&p^7&p^8\\
0&0&p^5&p^6&p^7\cr
\end{array}\right], \]
\[\left[\begin{array}{ccc|cc}
0&0&0&0&0\\
0&0&0&0&p^9\\
0&0&0&p^7&p^9\\
\hline
0&0&p^5&p^6&p^7\\
0&p^4&p^5&p^5&p^6\cr
\end{array}\right], \left[\begin{array}{ccc|cc}
0&0&0&0&p^9\\
0&0&0&p^7&p^9\\
0&0&p^5&p^7&p^8\\
\hline
0&p^4&p^5&p^5&p^6\\
p^2&p^4&p^4&p^4&p^5\cr
\end{array}\right], \
\left[\begin{array}{ccc|cc}
0&0&0&p^7&p^9\\
0&0&p^5&p^7&p^8\\
0&p^4&p^5&p^6&p^8\\
\hline
p^2&p^4&p^4&p^5&p^6\\
p^2&p^3&p^4&p^4&p^4\cr
\end{array}\right], \]
\[\left[\begin{array}{ccc|cc}
0&0&p^5&p^7&p^8\\
0&p^4&p^5&p^6&p^8\\
p^2&p^4&p^4&p^6&p^7\\
\hline
p^2&p^3&p^4&p^4&p^6\\
p&p^3&p^3&p^4&p^4\cr
\end{array}\right],
\left[\begin{array}{ccc|cc}
0&p^4&p^5&p^6&p^8\\
p^2&p^4&p^4&p^6&p^7\\
p^2&p^3&p^4&p^5&p^7\\
\hline
p&p^3&p^3&p^4&p^5\\
p&p^2&p^3&p^3&p^3\cr
\end{array}\right], \left[\begin{array}{ccc|cc}
p^2&p^4&p^4&p^6&p^7\\
p^2&p^3&p^4&p^5&p^7\\
p&p^3&p^3&p^5&p^6\\
\hline
p&p^2&p^3&p^3&p^5\\
*&p^2&p^2&p^3&p^3\cr
\end{array}\right]\]
\[ \left[\begin{array}{ccc|cc}
p^2&p^3&p^4&p^5&p^7\\
p&p^3&p^3&p^5&p^6\\
p&p^2&p^3&p^4&p^6\\
\hline
*&p^2&p^2&p^3&p^4\\
*&p&p^2&p^2&p^2\cr
\end{array}\right], \left[\begin{array}{ccc|cc}
p&p^3&p^3&p^5&p^7\\
p&p^2&p^3&p^4&p^6\\
*&p^2&p^2&p^4&p^5\\
\hline
*&p&p^2&p^2&p^4\\
*&p&p&p^2&p^2\cr
\end{array}\right], \left[\begin{array}{ccc|cc}
p&p^2&p^3&p^5&p^7\\
*&p^2&p^2&p^4&p^5\\
*&p&p^2&p^3&p^5\\
\hline
*&p&p&p^2&p^3\\
*&*&p&p&p\cr
\end{array}\right]\]
\[\left[\begin{array}{ccc|cc}
p&p^2&p^3&p^5&p^7\\
*&p&\overbrace{p^2}&p^3&p^5\\
*&\overbrace{p}&p&\overbrace{p^3}&p^4\\
\hline
*&*&\overbrace{p}&p&\overbrace{p^3}\\
*&*&*&\overbrace{p}&p\cr
\end{array}\right],
\mathcal J_{14}=\left[\begin{array}{ccc|cc}
p&p^2&p^3&p^5&p^7\\
*&p&p&p^3&p^5\\
*&*&p&p^2&p^4\\
\hline
*&*&*&p&p^2\\
*&*&*&*&p\cr
\end{array}\right]=\J\]

 with annihilating length $d=14$.
\end{example}
\section{The (uas)-(ucs)bounded  abelian  $p$--groups} \label{nil}
Let $\{J_t\}$ be the upper annihilating sequence of $\J$.  We denote
the nilpotency class of $\Delta(G) $ by $n(\Delta)$, and $ Z
=\Delta(G) \cap \mathcal Z$, where  $\mathcal Z$ is the center of
the ring $\E$, (\ref{Prel}.b). Denoting by $\Gamma
_t=(\J_t+1)\mathcal Z$, for $t>1$, and $\Gamma _0=\{1\}$, we will
prove that $\Gamma _t=\mathcal Z_t$. By Theorem \ref{uas=ucs} the
series $(\Gamma_t)$ is a central series, so we only need to prove
that $\mathcal Z _t\leq \Gamma _t$, for all $t$. The  Lemma
\ref{chac} gives a characterization  of the elements in the groups
$\Gamma _t$. This characterization of $\Gamma _t$ is very important
because we will prove that $\Gamma _t=\mathcal Z_t$ for almost all
of groups $G$.
\begin{definition} We will say that the abelian  p-group $G$ has the property
(uas)-(ucs)--group for the Jacobson radical $\J$ of $\E$, if the
upper central series $(\mathcal Z_t)_{t=0,n(\Delta)}$ of the
$p$-subgroup $\Delta(G)$ satisfies the following $\mathcal Z
_t=\Gamma_t$ for all $t$ from $1$ to  $n(\Delta)$.
\end{definition}

If $1_{rr} $ is the identity matrix of size $r\times r$, we use the
following notation for the elements of $\Gamma _t$: $C+1_{rr}\in
\Gamma _t$, then $C+1_{rr}=(c^{(l,k)}_{ij})_{r\times r}+1_{rr}$. On
the other hand $\Gamma _t=(\J_t+1)\mathcal Z$.

\begin{lemma}[Characterization of $\Gamma_t$]\label{chac}  \hspace{2.6in}

\emph{( 1 )} $C+1_{rr}\in \Gamma _t$ if, and only if
$C-c_{ss}^{(r,r)} 1_{rr}\in J_t$,
\[\emph{( 2 )}\hspace{.4in} \Gamma _t=\left\{ C+1_{rr}\in \Delta (G)| \left\{\begin{array}{ll}
c_{ij}^{(l,k)}\equiv 0 \pmod{p^{n_j-\alpha{(i,j,t)}}}& i\ne j;\\
c_{ii}^{(l,l)}\equiv c_{ss}^{(r,r)} \pmod {p^{n_i-\alpha{(i,i,t)}}}& \\
c_{ii}^{(l,k)}\equiv 0 \pmod{p^{n_i-\alpha{(i,i,t)}}}& l\ne k\cr
\end{array} \right\}\right\}  \]
\end{lemma}
\begin{proof}
(1) We know that $C+1_{rr}\in \Gamma _t=\J_t+Z_\Delta$ if, and only
if there exists $\overline C\in \J_t$, and a scalar matrix
$(pc+1)1_{rr}\in  Z$ such that $\overline C+(pc+1)1_{rr}=C+1_{rr}$.
 If $C+1_{rr}\in \Gamma _t$, then either for $i\ne j$, or $i=j$,
and $l\ne k$, we have that $c^{(l,k)}_{ij}=\overline
c^{(l,k)}_{ij}\equiv 0 \pmod {p^{n_j-\alpha{(i,j,t)}}}$. Because
$c_{jj}^{(l,l)}=\overline c_{jj}^{(l,l)}+pc$, for all $l$, and $j$.
We have that
\[c_{jj}^{(l,l)}- c_{ss}^{(r,r)}=\overline c_{jj}^{(l,l)}- \overline
c_{ss}^{(r,r)}\equiv 0 \pmod {p^{n_j-\alpha{(j,j,t)}}},\] by
Proposition \ref{prop:uf}, property (4).

 Therefore
$C-c_{ss}^{(r,r)}1_{rr}\in \J_t$, and
$C+1_{rr}=(C-c_{ss}^{(r,r)}1_{rr})+ (c_{ss}^{(r,r)}+1)1_{rr})\in
\Gamma _t$, if and only if $C-c_{ss}^{(r,r)}1_{rr}\in \J_t$

(2) It is a consequence of part (1).
\end{proof}

\section{Construction of the upper central series using the upper annihilating sequence}\label{algorithm}
This section is devoted to prove  Theorem \ref{p>2}.  In fact, we
want to prove  $\Gamma _t=\mathcal Z_t$ for all $t\geq 1$, we know
that $\Gamma _t\leq \mathcal Z_t$, so we need to prove that
$\mathcal Z_t \leq \Gamma _t$. We  use induction on $t$, for $t=1$
is proved in Lemma \ref{Center}, under assumption  the property
holds for $t-1$, then  we  prove the property for $t$ holds too. It
is done in
 Lemmas \ref{cond:CR}, and \ref{cases}.  We want to remark that in
 the hypothesis of induction on $t$  is $\Gamma
_{t-1}=\mathcal Z_{t-1}$.
\begin{theorem}\label{p>2}
The group $G=\oplus_{i=1}^s G_i$, where the $p$-rank of $G_i$ is
$r_i$, is a (uas)-(ucs)--group for $\Delta (G)$, in the following
cases
\begin{itemize}
\item[(1)] $r_s>1$, and $s\geq 2$, for all prime number $p$,
\item [(2)] $\sigma (G)\geq 2$, for all prime number $p$,
\item [(3)]  $p\geq 3$.
\end{itemize}
\end{theorem}

\begin{lemma}\label{Center}
The center $\mathcal Z_1$ of the subgroup $\Delta (G)$ is equal
$\Gamma _1=(\J_1+1)\mathcal Z$.
\end{lemma}
\begin{proof}
We know that $\Gamma _1\leq  \mathcal Z_1$ by Theorem \ref{uas=ucs}.
So, we will prove  $\mathcal Z_1\leq \Gamma _1$. Taking two elements
in $\Delta (G)$ and considering $A+1\in \mathcal Z_1$, we have that
\begin{equation} \begin{array}{l}
(A+1)(B+1)=(B+1)(A+1)\\
BA-AB=(A,B)=0,\\
\end{array}
\end{equation} for all $B+1\in \Delta (G)=\J+1$. But
\begin{equation}\label {C1}(BA-AB)_{ij}=\sum_{k=1}^s(B_{ik}A_{kj}-A_{ik}B_{kj})\equiv 0 \pmod{p^{n_j}}.\end{equation} Taking for
the matrix $B=(B_{ij})_{s\times s}$ the elementary matrices
$A-a_{ss}^{(r,r)}1_{rr}\in \J_1$  so  $A+1\in  \Gamma _1$. Therefore
$\Gamma _1= \mathcal Z_1$, and our claim holds.
\end{proof}
A trivial consequence is that $\Delta (G)$ is abelian if, and only
if $G$ has type $t(G)=(2).$

 In the following Lemma we
prove that $\mathcal Z_t\leq \Gamma _t$, for the cases described in
the Lemma \ref{cases}, assuming the condition (\ref{CR}).
\begin{lemma}\label{cond:CR}
The group $G$ is a (uas)-(ucs)--group for $\Delta (G)$ if and only
if the following  condition  holds,
\begin{equation}\label{CR}
If  \ A+1\in \mathcal Z_t \ then \ BA-AB=(B,A)\in \J_{t-1} \ for \
all \ B+1\in \Delta (G).\end{equation}
\end{lemma}
\begin{proof} ($\Rightarrow$ ) Suppose the group is (uas)-(ucs)-group, that is $\Gamma _t=\mathcal Z _t$ for all $t$, then the
condition holds by Corollary \ref{COR}.

($\Leftarrow$) Suppose the group satisfies the condition
$(\ref{CR})$.
 By definition of $\mathcal Z_t$,
we have that for all $B+1\in \Delta(G)$, and $A+1\in \mathcal Z_t$
then \[[B+1,A+1]=C+1\in \mathcal Z_{t-1}=\Gamma
_{t-1}=(\J_{t-1}+1)Z=\J_{t-1}+\mathcal Z\] by induction. This
implies that $BA-AB=C(A+1)(B+1)$, where $C+1\in \mathcal
Z_{t-1}=\Gamma _{t-1}$. Therefore
\[(C(A+1)(B+1))_{ij}=\sum_{k=1}^s
C_{ik}\sum_{l=1}^s(A_{kl}+\delta_{kl})(B_{lj}+\delta_{lj}),
\]
But $(BA-AB)_{kj}=\sum_{i=1}^s(B_{ki}A_{ij}-A_{ki}B_{ij})\equiv 0
\pmod{ \J_{t-1}}$, by condition $(\ref{CR})$. So,
\begin{equation}\label{equ1}
\hspace{-.6in}{(kj)}\hspace{.5in}\sum_{i=1}^s(B_{ki}A_{ij}-A_{ki}B_{ij})\equiv
0 \pmod {p^{n_j-\alpha (k,j,t-1)}}
\end{equation}
Now, we will prove that $A-a_{ss}^{(r,r)}1_{rr}\in \J_t$, and of
course $A+1\in \Gamma _t$ then the  group is (uas)-(ucs)--group. The
idea  is to select for the element $B$  different matrices, in such
a way  that all the conditions for the upper annihilating functions
are satisfied. In fact, in (\ref{equ1})-(kj),we have the following
\[\begin{array}{|ccc|}
\hline
 taking\ B=E_{ki}^{(u,v)} & for\ i\ne j, &
e_{ki}^{(u,v)}a_{ij}^{(v,l)}\equiv 0 \pmod{p^{n_j-\alpha
(k,j,t-1)}}\\
if \ k=i+1 \ then & e_{i+1,i}^{(u,v)}=1, & a_{ij}^{(v,l)}\equiv 0
\pmod{p^{n_j-\alpha (i+1,j,t-1)}}\\
if \ k=i-1 \ then & e_{i-1,i}^{(u,v)}=p^{n_i-n_{i-1}}, &
a_{ij}^{(v,l)}\equiv 0
\pmod{p^{n_j-n_i+n_{i-1}-\alpha (i-1,j,t-1)}}\\
\hline taking\ B=E_{kj}^{(u,v)} & for\ i\ne j, &
a_{ki}^{(l,u)}e_{i,j}^{(u,v)}\equiv 0 \pmod{p^{n_j-\alpha
(k,j,t-1)}}\\
if \ j=i+1 \ then & e_{i,i+1}^{(u,v)}=p^{n_{i+1}-n_i}, &
a_{ki}^{(l,u)}\equiv 0
\pmod{p^{n_{i+1}-n_{i+1}+n_i-\alpha (k,i+1,t-1)}}\\
if \ j=i-1 \ then & e_{i,i-1}^{(u,v)}=1, & a_{ki}^{(l,u)}\equiv 0
\pmod{p^{n_{i-1}-\alpha (k,i-1,t-1)}}\\
\hline
\end{array}
\]
Similarly, for $i\ne j$ we obtain  that $a_{ij}^{(l,k)}\equiv 0
\pmod {p^{n_j-1-\alpha (i,j,t-1)}}$, and $a_{ij}^{(l,k)}\equiv 0
\pmod {p^{n_j-\alpha (i,j,t)}}$ because
\[\begin{array}{|l|}
\hline
 \alpha{(i,j,t)}=\bigwedge
\{n_j-n_{j-1}+\alpha{(i,j-1,t-1)};n_i-n_{i-1}+\alpha{(i-1,j,t-1)}; \\
\alpha{(i,j,t-1)}+1;\alpha{(i+1,j,t-1)};\alpha{(i,j+1,t-1)}\},\cr
\hline
\end{array}\]
In the same way for the matrix $A-a_{ss}^{(r,r)}1_{rr}$ we prove
that
\[a_{jj}^{(l,l)}-a_{ss}^{(r,r)}\equiv 0 \pmod {p^{n_j-\alpha (j,j,t)}},\ for \ j<s\] and
if the $p$-rank $r_s$ of the last homocyclic  group $G_s$ is greater
than $1$, we have that \[a_{ss}^{(l,l)}-a_{ss}^{(r,r)}\equiv 0 \pmod
{p^{n_s-\alpha (s,s,t)}}, \ for \ l\leq r .\] We want to remark that
in the place $rr$ of the matrix $A-a_{ss}^{(r,r)}1_{rr}$ we have
$0$. Therefore $A-a_{ss}^{(r,r)}1_{rr}\in\J_t$, $A+1\in \J_t+1$, and
the group is a (uas)-(ucs)-group. So, our claim holds.
\end{proof}
\begin{lemma}\label{cases}
The condition (\ref{CR}) holds in the following cases
\begin{itemize}
\item[(1)] $r_s>1$, and $s\geq 2$, for all prime number $p$,
\item [(2)] $\sigma (G)\geq 2$, for all prime number $p$,
\item [(3)] otherwise, only  for $p\geq 3$.
\end{itemize}
\end{lemma}
\begin{proof}
I need to prove that  $(B,A)\in \J _{t-1}$, for  $A+1\in \mathcal Z
_t$, and $B+1\in \Delta (G)$. We have the following expression
\[(B,A)=C(A+1)(B+1)=\sum_{k=1}^s
C_{ik}(\sum_{l=1}^s(A_{kl}+\delta_{kl})(B_{lj}+\delta_{lj})),
\]
we have that $(B,A)\in \J_{t-1}$ if, and only if $C\in \J_{t-1}$. We
know that for all $B$ and $A+1\in \mathcal Z_t$, we have $C+1\in
\mathcal Z _{t-1}=\Gamma _{t-1}$, by induction. So by Lemma
\ref{chac}, we have that: $C+1\in \Gamma _{t-1} \ \Leftrightarrow \
C- c_{ss}^{(r,r)}1_{rr}\in \J _{t-1}$, and  $C\in \J_{t-1}$ if and
only if $c_{ss}^{(r,r)}1_{rr}\in \J _{t-1}$. But
\[BA-AB=(C-c_{ss}^{(r,r)}1_{rr}
+c_{ss}^{(r,r)}1_{rr})(A+1)(B+1).\] By $C+1\in \mathcal Z
_{t-1}=\Gamma _{t-1}$, then $C- c_{ss}^{(r,r)}1_{rr}\in \J _{t-1}$,
so
\[BA-AB\equiv c_{ss}^{(r,r)}1_{rr}(A+1)(B+1) \pmod {\J_{t-1}}.\]

Because the elementary endomorphism defined in Section \ref{Prel}
form a set of generator for the Jacobson radical $\J$, we will prove
that our claim is true for all the elementary matrix in $\J$.
\[(E_{ij}^{(k,u)},A)\equiv c_{ss}^{(r,r)}\sum_{l=1}^s
(A_{il}+\delta_{sl})(E^{(k,u)}_{l,j}+\delta_{lj}), \pmod
{p^{n_j-\alpha(i,j,t-1)}}.
\]

We want to remark here that in general,  all elementary endomorphism
have a different matrix $C+1$. The notation of the entries in the
matrix $C$ are the following $c_{i,j}^{(l,k)}$, so if we need to
include another notation in these elements, it will be more
complicated. Then we prefer to consider only the notation
$c_{ss}^{(r,r)}$ for all the matrices, knowing  that in each case we
have the possibility to have a different $C$. We will prove that in
all the cases the element $c_{ss}^{(r,r)}\equiv 0 \pmod
{p^{n_s-\alpha (s,s,t-1)}}$ First we prove that for all  the
elementary matrices  such that $(i,j)\ne (s,s)$ then the entry
$c_{s,s}^{(r,r)}$ in the associated matrix $C$ satisfies the
condition for the elements in $\J_{t-1}$, for $k\ne r,$ and $u\ne r$
\[\begin{array}{|c|c|}
\hline    (E_{ij}^{(k,u)},A) _{ss}^{(r,r)}=0& then \  0\equiv
c_{s,s}^{(r,r)}(a_{ss}^{(r,r)}+1)
\pmod {p^{n_s-\alpha(s,s,t-1)}} \\
\hline but \  (a_{ss}^{(r,r)}+1)^{-1} \ exists &  then \ 0\equiv
c_{s,s}^{(r,r)}\pmod {p^{n_s-\alpha(s,s,t-1)}}\\
\hline then \ c_{s,s}^{(r,r)}1_{rr}\in  \J_{t-1} &  and \ (E_{ij}^{(k,u)},A)\in \J_{t-1}\\
 \hline
\end{array}\]
Then we have proved the condition for all the entries in the
elementary endomorphism, where $k<r$, and $u<r$.

 \textbf{(Case number 1)} $r_s>1$, and $s\geq 2$, for all prime number
 $p$.

 Suppose that $k=r$ and $l\ne r$. Then in the place
 $(r-1,r-1)$ we have $0\equiv
c_{s,s}^{(r,r)}(a_{ss}^{(r-1,r-1)}+1)\pmod
{p^{n_s-\alpha(s,s,t-1)}}$, and our claim holds because
$(E_{sj}^{(r,u)},A)\in \J_{t-1}$. Similarly we prove that
$(E_{is}^{(l,r)},A)\in \J_{t-1}$, even $(E_{ss}^{(r,r)},A)\in
\J_{t-1}$. The case is solved.

\textbf{(Case number 2)} $\sigma (G)\geq 2$, for all prime number
$p$.

The solution is similar the above case, considering $r_s=1, $ and
$\sigma (G)\geq 2$, that are the correspondent unsolve case.

\textbf{(Case number 3)} $p\geq 3$.

 We have the cases for giving solution\\
(a) $r_s=1$, $\sigma (G)=1$, and $n_s>n_{s-1}+1$. The solution is
similar the first case,  \\
(b) $r_s=1$, $\sigma (G)=1$, and $n_s=n_{s-1}+1$, we need the
condition $p\geq 3$, and only in this case.

We give here the solution:
\[\begin{array}{|c|}
\hline
(E_{sj}^{(r,u)},A)_{ss}^{(r,r)}=e_{sj}^{(r,u)}a_{js}^{(u,r)}\ \Rightarrow \\
 e_{sj}^{(r,u)}a_{js}^{(u,r)}\equiv
c_{s,s}^{(r,r)}(a_{ss}^{r,r}+1)
\pmod {p^{n_s-\alpha(s,s,t-1)}}\ (I) \\
\hline (E_{sj}^{(r_s,u)},A)
_{jj}^{(r,r)}=-a_{js}^{(u,r)}e_{sj}^{(r,u)} \ \Rightarrow \\
-a_{js}^{(u,r)}e_{sj}^{(r,u)}\equiv
c_{s,s}^{(r,r)}(a_{jj}^{(l,l)}+1+ a_{js}^{(u,v)})
\pmod {p^{n_{j}-\alpha(j,j,t-1)}},\ (II)\\
\hline
taking \ j=s-1 \ we \ have \ by Lemma \ref{alpha}\\   n_s-\alpha(s,s,t-1)=n_{s-1}-\alpha (s-1,s-1,t-1) \\
 \hline
 summand \ (I) \ and \ (II) \\
  0\equiv c_{s,s}^{(r,r)}(a_{jj}^{(l,l)}+1+ a_{js}^{(u,v)}+a_{ss}^{(r,r)}+1) \pmod
 {p^{n_s-\alpha(s,s,t-1)}}\\
 \hline
because\ (a_{jj}^{(l,l)}+ a_{js}^{(u,v)}+a_{ss}^{(r,r)}+2)^{-1}\  exists \   for \  p\ne 3 \\ c_{s,s}^{(r,r)}1_{ss}\in \J_{t-1}, \ and  \ (E_{sj}^{(r,u)},A)\in \J_{t-1} \\
 \hline
\end{array}\]

Similarly we prove that $(E_{is}^{(l,r)},A)\in \J_{t-1}$, even
$(E_{ss}^{(r,r)},A)\in \J_{t-1}$. The case is solved.

\end{proof}

 If $S=\{1,\ldots,s\}$, we
denote by $y(G)$ the least $t$ such that $\alpha{(i,j,t)}\ge
f_\J(i,j)$, for $(i,j)\in S \times S \setminus \{(s,s)\}$. Clearly
$y(G)\leq l(\J)=d$
\begin{theorem}\label{stop}
If $G$ is a (uas)-(ucs)--group for $\Delta(G)$ then
\begin{equation}n(G)=\left\{\begin{array}{ccc} l(\J) & if  & r_s\geq 2\\
y(G)& if & r_s=1\cr
\end{array} \right\}
\end{equation}
\end{theorem}
\begin{proof}
By Theorem \ref{uas=ucs}, we know  that $n(G)\leq l(\J)$.

Suppose $r_s>2$, and $n(G)=d=l(\J)$ then we have that $\alpha
{(s,s,d-1)}<n_s-1$.

 By definition of  elementary endomorphism   the matrix $E^{(l,k)}_{ss}\in \J_d=\J$
  for $(l,k)\ne (r,r)$. But $E^{(l,k)}_{ss}\not \in \J_{d-1}$
   because $E^{(l,k)}_{ss}\not \equiv 0\pmod {p^{n_s-\alpha{(s,s,d-1)}}}$, because  $n_s-\alpha{(s,s,d-1)}\ge
   2$. So, $E^{(l,k)}_{ss}+1_{ss} \in \Delta(G)  - \Gamma_{d-1},$ and as a consequence $\Gamma
   _{d-1}< \Delta$, and  $n(G)=d=l(\J)$.

   If $r_s=1$, we prove that $\Delta (G) =\J+1 \leq \Gamma
   _{y(G)}$, because trivially $\Gamma _{y(G)}\leq \Delta(G)$. If we consider $C+1_{ss}\in \Delta(G)=\Gamma _d$,   then $C-c_{ss}^{(r,r)}1_{rr}\in \J_d$,
   but we assume that $C-c_{ss}^{(r,r)}1_{rr}\in \J_{y(G)}$. Therefore $C+1_{ss}\in
   \Gamma _{y(G)}$ by Lemma \ref{chac}, then $\Delta (G) \leq \Gamma_{y(G)}$, and our claim holds.
\end{proof}

\section{Algorithm for computing the upper central series , and the nilpotency class of $\Delta
(G)$}\label{Algorithm1}Let $G$ be a group with the property
(uas)-(ucs)-group, that is the hypercenters $\mathcal Z_t=\Gamma
_t=(A_t+1)Z$, where $\{\J_t\}$ is the upper annihilating sequence of
$\J$, and $Z=\Delta(G)\cap \mathcal Z$. Let $C+1\in \mathcal
Z_t=\Gamma _t $ be the matrix described as follows:
$C=(c_{ij}^{(l,k)})_{r\times r}$, with the entries in $\Z_{p^{n_j}}$
for $1\leq j\leq s$.
\begin{algorithm} \label{algorithm} Input: $\underline t(G)=(n_1,\ldots,
n_s)$, and the $p$-ranks  ($r_1$, $\ldots$, $r_s$).
\begin{itemize}
\item[(A1)] Compute $r=r_1+\cdots +r_s$.
\item [(A1)]Write the matrix $C^t=(c_{ij}^{(l,k)})_{r\times r}$, with the entries
 described  above.
\item[(A2)]  If $r_s>1$  for $t=0$ until least $t$ such that $\alpha(i,j,t)\geq f_\J(i,j)$ for all
$(i,j)$ do
 \begin{itemize}
 \item[] else for $t=1$ until the least $t$ such that $\alpha(i,j,t)\geq f_\J(i,j)$ for all
$(i,j)\ne (s,s))$ do
\item[1.]  $\alpha (i,j,0)=0$ for all $(i,j)$.
\item[2.] Compute the function $\alpha(i,j,t)$ using formulas in
Corollary
 \ref{formulas}.
\end{itemize}
 \item[(A4)] Write on  the matrix  $C^t$ of $\mathcal Z_t$,  the following:
\[ \ for \ i\ne j \ do \ c_{i,j}^{(l,k)}=p^{n_j-\alpha(i,j,t-1)}c_{i,j}^{(l,k)}\rightarrow
p^{n_j-\alpha(i,j,t)}c_{i,j}^{(l,k)}\]
\[ \ for \ l\ne k\ do \ c_{j,j}^{(l,k)}=p^{n_j-\alpha(j,j,t-1)}c_{j,j}^{(l,k)}\rightarrow
p^{n_j-\alpha(j,j,t)}c_{j,j}^{(l,k)}\]
\item[(A5)] Write the relations for the elements in the diagonal for $j\leq s$, and $ l\ne
r$ do
\[ \ c_{jj}^{(l,l)}=p^{n_j-\alpha(j,j,t-1)}\overline c_{jj}^{(l,l)}+pc_{ss}^{(r,r)}\rightarrow
c_{jj}^{(l,l)}=p^{n_j-\alpha(j,j,t)}\overline
c_{jj}^{(l,l)}+pc_{ss}^{(r,r)} ,\]
  and   $pc_{ss}^{rr}\longrightarrow pc_{ss}^{rr}$,
\item[(A6)] Save $C^t$, and write $\mathcal Z_t=\{C^t+1_{rr}\}$.
\item[(A7)] Save and write the number  $t=n(G)$ when  the algorithm stop.
\end{itemize}
Output: $n(G)=$the nilpotency class of $\Delta (G)$,  and the
centers $\mathcal Z_1$, $\ldots$, $\mathcal Z_{n(G)}$.
\end{algorithm}
\begin{example} For a group of type
$(3,5,7)$, with ranks $(1,1,2)$, when $p\geq 3$, we have  $\Gamma
_t=\mathcal Z _t$. Applying  the Algorithm we found   the upper
central series of $\Delta (G)$ Then the elements of the upper
central series are
\[\mathcal Z_1=\left\{\left[\begin{array}{c|c|cc}
pc_{44}+1&0&0&0\\
\hline
0&pc_{44}+1&0&0\\
\hline
0&0&p^6c_{33}+pc_{44}+1&p^6c_{34}\\
0&0&p^6c_{43}& c_{44}+1\cr
\end{array}\right] \right\},\] \[\mathcal Z_2=\left\{\left[\begin{array}{c|c|cc}
pc_{44}+1&0&0&0\\
\hline
0&pc_{44}+1&p^6c_{23}&p^6c_{24}\\
\hline
0&p^4c_{32}&p^5c_{33}+pc_{44}+1&p^5c_{34}\\
0&p^4c_{42}&p^5c_{43}&pc_{44}+1\cr
\end{array}\right] \right\} \]
\[\mathcal Z_3=\left\{\left[\begin{array}{c|c|cc}
pc_{44}+1&0&p^6c_{13}&p^6c_{14}\\
\hline
0&p^4c_{22}+pc_{44}+1&p^5c_{23}&p^5c_{24}\\
\hline
 p^2c_{31}&p^3c_{32}&p^4c_{33}+pc_{44}+1&p^4c_{34}\\
 p^2c_{41}&p^3c_{42}&p^4c_{43}& pc_{44}+1\cr
\end{array}\right] \right\} \]
\[\mathcal Z_4=\left\{\left[\begin{array}{c|c|cc}
pc_{44}+1&p^4c_{12}&p^5c_{13}&p^5c_{14}\\
\hline
p^2c_{21}&p^3c_{22}+pc_{44}+1&p^4c_{23}&p^4c_{24}\\
\hline
pc_{31}&p^2c_{32}&p^3c_{33}+pc_{44}+1& p^3c_{34}\\
pc_{41}&p^2c_{42}& p^3c_{43}&pc_{44}+1 \cr
\end{array}\right]  \right\} \]
\[\mathcal Z_5=\left\{\left[\begin{array}{c|c|cc}
p^2c_{11}+pc_{44}+1&p^3c_{12}&p^4c_{13}&p^4c_{14}\\
\hline
pc_{21}&p^2c_{22}+pc_{44}+1&p^3c_{23}&p^3c_{24}\\
\hline
 c_{31}&pc_{32}&p^2c_{33}+pc_{44}+1&p^2c_{34}\\
 c_{41}&pc_{42}&p^2c_{43}&pc_{44}+1 \cr
\end{array}\right]  \right\} \]
\[\Delta (G)=\mathcal Z_6=\left\{\left[\begin{array}{c|c|cc}
pc_{11}+1&p^2c_{12}&p^4c_{13}&p^4c_{14}\\
\hline
c_{21}&pc_{22}+1&p^2c_{23}&p^2c_{34}\\
\hline
 c_{31}&c_{32}&pc_{33}+1&pc_{34}\\
 c_{41}&c_{42}&pc_{43}&pc_{44}+1 \cr
\end{array}\right] \right\}. \]
The class of nilpotency of the group $\Delta (G)$ is  $n(G)=6$,
\end{example}

\begin{example} For a group of type
$(2,4,7)$, with ranks $(1,1,1)$, when $p\geq 3$, we have  $\Gamma
_t=\mathcal Z _t$ for $t=1$ to $t=y(G)$. Applying  the Algorithm we
found the upper central series of $\Delta (G)$ Then the elements of
the upper central series are
\[\mathcal Z_1=\left\{\left[\begin{array}{ccc}
pc_{33}+1&0&0\\
0&pc_{33}+1&0\\
0&0&pc_{33}+1\\
\end{array}\right]\right\}\]
\[\mathcal Z_2=\left\{\left[\begin{array}{ccc}
pc_{33}+1&0&0\\
0&pc_{33}+1&p^6c_{23}\\
0&p^3c_{32}&pc_{33}+1\\
\end{array}\right]  \right\} \]
\[\mathcal Z_3=\left\{\left[\begin{array}{ccc}
pc_{33}+1&0&p^6c_{13}\\
p^2c_{21}&p^3c_{22}+pc_{33}+1&p^5c_{23}\\
pc_{31}&p^2c_{32}&pc_{33}+1\\
\end{array}\right] \right\} \]
\[\mathcal Z_4=\left\{\left[\begin{array}{ccc}
pc_{33}+1&p^3c_{12}&p^5c_{13}\\
pc_{21}&p^2c_{22}+pc_{33}+1&p^4c_{23}\\
c_{31}&pc_{32}&pc_{33}+1\\
\end{array}\right]  \right\} \]
\[\mathcal Z_5=\left\{\left[\begin{array}{ccc}
pc_{11}+1&p^2c_{12}&p^5c_{13}\\
c_{21}&pc_{22}+1&p^3c_{23}\\
c_{31}&c_{32}&pc_{33}+1\\
\end{array}\right] \right\}=\Delta (G) \]
The class of nilpotency of the group $\Delta (G)$ is  $n(G)=5$, we
can check that in this case the length of the upper annihilating
sequence is $l(\J_t)=6$
\end{example}
\section*{Acknowledgments}
I want to thank Professor Laci Kovacs, from Sidney, Australia, for
informing me that my solution of the nilpotency class of $\Delta(G)$
 was the first  and only solution given until now to this problem,
open for more than 70 years. I would also like to thank  Professors
John C. McConnell, and Ken Goodearl, Ring Theorists, for informing
me that the concept of annihilating sequences  was new in
mathematics in 1996. I want to thank Professor Raymundo Bautista for
his useful suggestion. on the solution, using the algorithm and
formulas for the upper central series and the nilpotency class of
$\Delta(G)$, worked students from Sidney, Mexico City, and Santa
Clara, Cuba, so I want to include them in my Acknowledgments.

\end{document}